\documentclass{article}%
\usepackage{amsfonts}
\usepackage{amsmath}
\usepackage{amssymb}
\usepackage{hyperref}
\usepackage{graphicx}%
\RequirePackage[a4paper,top=2.54cm,bottom=2.54cm,left=1.90cm,right=1.90cm,%
                headsep=1em,includehead,includefoot]{geometry}
\RequirePackage{float}
\RequirePackage{indentfirst}
\providecommand{\U}[1]{\protect\rule{.1in}{.1in}}
\newtheorem{theorem}{Theorem}

\newtheorem{definition}[theorem]{Definition}

\newtheorem{proposition}[theorem]{Proposition}
\newtheorem{remark}[theorem]{Remark}

\newcommand{\R}{\mathbb{R}}
\newenvironment{proof}[1][Proof]{\noindent\textbf{#1.} }{\ \rule{0.5em}{0.5em}}

\providecommand{\keywords}[1]
{
  \small	
  \textbf{\textit{Keywords:}} #1
}

\providecommand{\MSC}[2]
{
  \small	
  \textbf{\textit{MSC:}} #1
}

\makeatletter
\def\@makefnmark{}
\makeatother
\newcommand{\myfootnote}[2]{\footnote{\textbf{#1}: #2}}

\begin{document}
\title{On the Hausdorff Measure of $\R^n$ with the Euclidean Topology}
\date{}
\author{Marco Bagnara\myfootnote{Marco Bagnara}{Scuola Normale Superiore, Piazza dei Cavalieri 7, 56126 Pisa, Italy, e-mail: marco.bagnara@sns.it}, Luca Gennaioli\myfootnote{Corresponding Author: Luca Gennaioli}{SISSA, Via Bonomea, 265, 34136 Trieste, Italy, e-mail: luca.gennaioli@sissa.it}, Giacomo Maria Leccese \myfootnote{Giacomo Maria Leccese}{SISSA, Via Bonomea, 265, 34136 Trieste, Italy, e-mail: giacomomaria.leccese@sissa.it}, Eliseo Luongo \myfootnote{Eliseo Luongo}{Scuola Normale Superiore, Piazza dei Cavalieri 7, 56126 Pisa, Italy, e-mail: eliseo.luongo@sns.it}}
\maketitle 
\begin{abstract}
  In this paper we answer a question raised by David H. Fremlin about the Hausdorff measure of $\R^2$ with respect to a distance inducing the Euclidean topology. In particular we prove that the Hausdorff $n$-dimensional measure of $\R^n$ is never $0$ when considering a distance inducing the Euclidean topology. Finally, we show via counterexamples that the previous result does not hold in general if we remove the assumption on the topology.
\end{abstract}
\\
\MSC{28A75, 28A78}

\section{Introduction}
The aim of this paper is to answer an open question stated by D. H. Fremlin in his famous book \cite{fremlin2000measure}, which can be found moreover on\\ \href{https://www1.essex.ac.uk/maths/people/fremlin/answer.pdf}{https://www1.essex.ac.uk/maths/people/fremlin/answer.pdf}:
\begin{equation}\label{long2}
\begin{split}
&\textit{let us consider a metric $\rho$ on $\mathbb{R}^2$ inducing the Euclidean topology,}\notag\\ &\textit{is it possible that $
    \mathcal{H}^2_{\rho}\left(\mathbb{R}^2\right)=0$?} 
\end{split}\tag{Q}
\end{equation}
By $\mathcal{H}^n_{\rho}$ we denote the n-dimensional Hausdorff measure according to Definition \ref{Hausdorff measure} below.
We give an answer to this problem in full generality, since our proof is valid in $\mathbb{R}^n,\ \forall n\geq 1$, showing that such a behaviour cannot happen. On the other hand, we will show in Remark \ref{pathologies} that, when the metric does not induce the usual Euclidean topology, counterexamples can be found. \\
Before stating our main theorem  in Section \ref{proof of thm}, we recall in this introductory section some classical tools for convenience of the reader (see \cite{federer2014geometric}, \cite{jost2005postmodern} for further details).

\begin{definition}[Hausdorff measure]\label{Hausdorff measure}
Let $(X,\textrm{d})$ be a metric space. We define the $n$-dimensional Hausdorff outer measure of $A\in\mathcal{P}(X)$ as
\begin{align}
    \mathcal{H}^n_{\text{d}}(A)&:=\operatorname{sup}_{\delta> 0}\mathcal{H}^n_{\delta,\text{d}}(A),\qquad\text{with} \\
    \mathcal{H}^n_{\delta,\text{d}}(A)&:=\inf\bigl\{\sum_{i\in I}\operatorname{diam}(A_i)^n\;:\;A\subseteq\cup_{i\in I}A_i,\;\operatorname{diam}(A_i)\leq\delta\bigr\},
\end{align}
where $\textrm{diam}(U)=\sup_{x,y\in U}\textrm{d}(x,y)$ and $I$ is an at most countable collection of indices.
\end{definition}
\begin{remark}
The usual definition of Hausdorff measure is given scaling the result by a dimensional constant that, for instance, in the Euclidean case is equal to $2^{-n}\omega _{n}$, where $\omega_{n}$ is the volume of the unit n-ball. We opted to overlook the constant in order to simplify the notation. Clearly Theorem \ref{main th} is not affected by this choice.
\end{remark}
To prove our result we will exploit the following well-known theorem.
\begin{theorem}[Dini]\label{dini}
Let $(K,d)$ be a compact metric space. Let $f_n:K\rightarrow \mathbb R$ be continuous functions such that \begin{align}
    f_n\leq f_{n+1}\qquad \forall n\in \mathbb{N}
\end{align}
and assume that \begin{align}
    f(x)=\lim_{n\rightarrow +\infty}f_n(x) \qquad \forall x\in K,
\end{align}
exists and the function $f:K\rightarrow\mathbb R$ is also continuous. Then $\left(f_n\right)_{n\in\mathbb{N}}$ converges uniformly to $f$ on $K$.
\end{theorem}
Moreover, we briefly recall the definition and some properties of the Brouwer Degree. See for instance \cite{dinca2021brouwer} for a complete treatment of this topic.
\begin{theorem}[Brouwer Degree]
There exists only one function, called Brouwer Degree and denoted by $\operatorname{deg}$, from the set of couples $(D, f )$, where
$D \subset \mathbb R^n$ is open and bounded and $f : \Bar{D} \to  \mathbb R^n$ is continuous with $0 \notin f (\partial D)$, into the set
$\mathbb{Z}$, which satisfies the following three properties:
\begin{itemize}
\item(Normalization) $\operatorname{deg}[\operatorname{id},D] = 1$ if $0 \in D$.
\item(Additivity) $\operatorname{deg}[f, D] = \operatorname{deg}[f, D_1] + \operatorname{deg}[f, D_2]$ if $D_1$ and $D_2$ are disjoint open
subsets of $D$ such that $0 \notin f (\Bar{D} \setminus (D_1 \cup D_2))$.
\item(Homotopy invariance) If $F \in C([0, 1]\times \Bar{D},  \mathbb R^n)$ and $0 \notin F ( [0, 1]\times\partial D)$,
then $\operatorname{deg}[F (t,\cdot), D]$ is independent of $t\in[0, 1]$.
\end{itemize}
\end{theorem}
\begin{definition}
If $D \subset  \mathbb R^n$ is open and bounded, $f \in C(\Bar{D},  \mathbb R^n)$ and $z \notin f (\partial D)$,
the Brouwer degree $\operatorname{deg}[f, D, z]$ is defined by
$\operatorname{deg}[f, D, z] = \operatorname{deg}[f (\cdot) - z, D]$.
\end{definition}
\begin{proposition}\label{prop degree}
If $z \notin f (\Bar{D})$, then $\operatorname{deg}[f, D, z] = 0.$\\ Equivalently, if $\operatorname{deg}[f, D, z] \not= 0$,
there exists at least one $x \in D$ such that $f (x) = z$.
\end{proposition}
\section{Main results}\label{proof of thm}We are now in the position to state our main theorem.
\begin{theorem}\label{main th}
Let $( \mathbb R^n,\rho)$ be a metric space with $\rho$ inducing the Euclidean topology, then $\mathcal{H}^n_\rho( \mathbb R^n)>0$.
\end{theorem}
\begin{proof}
Assume by contradiction that there exists a distance $\rho$ in $ \mathbb R^n$ such that $\mathcal{H}^n_\rho(\mathbb R^n)=0$. We denote by $\mathbb{B}(0,1)$ the closed unit ball with respect to Euclidean metric and we consider the identity map \begin{equation}
    \operatorname{id}:(\mathbb{B}(0,1),\rho)\longrightarrow(\mathbb{B}(0,1),\textrm{d}_{\textrm{eucl}}).
\end{equation}
Such a map is an homeomorphism by assumption, but it carries no metric information a priori. Let us write \begin{equation}
    \operatorname{id}(x)=\bigl(\pi_1(x),\dots,\pi_n(x)\bigr)
\end{equation} and define 
\begin{equation}
    \label{Lipschitz approx}
    \pi_i^\varepsilon(x):=\min_{z\in \mathbb{B}(0,1)}\biggl[\pi_i(z)+\frac{1}{\varepsilon}\rho(x,z)\biggr]\qquad\forall i =1,...,n\quad\forall x\in\mathbb{B}(0,1),
\end{equation}
where we are using that $\mathbb{B}(0,1)$ is compact also for the metric $\rho$. The latter functions are Lipschitz, since they are the infimum of a family of equi-Lipschitz functions, more precisely 
\begin{equation}
    |\pi_i^\varepsilon(x)-\pi_i^\varepsilon(y)|\leq\frac{1}{\varepsilon}\rho(x,y)\quad\forall x,y\in\mathbb{B}(0,1).
\end{equation}
Such functions converge pointwise to the components of the identity in the compact ball $\mathbb{B}(0,1)$ as $\varepsilon\to 0$. Indeed, consider a sequence $(z_\varepsilon)_{\varepsilon>0}\subseteq\mathbb{B}(0,1)$ such that \begin{equation}\label{equa:zeta_eps}
    \pi_i^\varepsilon(x)=\pi_i(z_\varepsilon)+\frac{1}{\varepsilon}\rho(x,z_\varepsilon).
\end{equation}
This sequence is bounded and by compactness it admits a convergent subsequence. Due to equation \eqref{equa:zeta_eps} and the bound 
\begin{equation}
    1\ge\pi_i\ge\pi_i^\varepsilon\ge-1,
\end{equation}
it follows that $\lim_{\varepsilon\to 0}\rho(z_\varepsilon,x)=0$, which means that the whole sequence converges to $x$, leading to the desired pointwise convergence. Now, since we have $\pi_i^\varepsilon(x)\geq \pi_i^{\varepsilon+\gamma}(x)$ for every $\gamma,\varepsilon>0$ and $\forall x\in\mathbb{B}(0,1)$, by Dini's theorem $\pi_i^\varepsilon$ converges uniformly to $\pi_i$ on $\mathbb{B}(0,1)$ for every $i=1,...,n$.
Summing up we have obtained a sequence 
\begin{equation}
    F^\varepsilon = (\pi_1^\varepsilon,...,\pi_n^\varepsilon): (\mathbb{B}(0,1),\rho)\longrightarrow (\mathbb R^n,\textrm{d}_{\textrm{eucl}})
\end{equation} such that 
\begin{equation}
    \textrm{d}_{\textrm{eucl}}\bigl(F^\varepsilon(x), F^\varepsilon(y)\bigr)\leq C_\varepsilon\rho(x,y)\quad\forall x,y\in\mathbb{B}(0,1)
\end{equation}
with $C_\varepsilon>0$ and such that it converges uniformly to the identity in $\mathbb{B}(0,1)$. The following claim is of crucial importance.\\ \emph{Claim}: there exists $\varepsilon>0$ such that $F^\varepsilon(\mathbb{B}(0,1))$ has non-empty interior.\\
Fix $\hat \varepsilon>0$ such that 
\begin{equation}\label{equa:bound}
    \sup_{x\in\mathbb B(0,1)}\textrm{d}_{\textrm{eucl}}(F^{\varepsilon}(x),x)\le \frac{1}{2}
\end{equation} for every $\varepsilon\in [0,\hat \varepsilon]$ and consider the function
\begin{equation}
    F:[0,\hat \varepsilon]\times \mathbb{B}(0,1)\to  \mathbb R^n
\end{equation} defined by the relation $F(\varepsilon,\cdot)=F^\varepsilon$ for $\varepsilon>0$ and $F(0,\cdot)=id$. We prove that the function $F$ is a continuous function, or in other words that $F$ is an homotopy between id and $F^{\hat\varepsilon}$.
First we observe that for every $\varepsilon_n\nearrow\varepsilon$ in $(0,\hat \varepsilon]$, given $z$ such that 
\begin{equation}
    F_i^{\varepsilon}(x)=\pi_i(z)+\frac{1}{\varepsilon}\rho(x,z),
\end{equation}
then
\begin{equation}
    \pi_i(z)+\frac{1}{\varepsilon_n}\rho(x,z)\ge F_i^{\varepsilon_n}(x)\ge F_i^{\varepsilon}(x)
\end{equation}
and taking the limit for $n\to+\infty$, we obtain that $\lim_{n\to+\infty}F_i^{\varepsilon_n}(x)=F_i^{\varepsilon}(x)$.
Also, for every $\varepsilon_n\searrow\varepsilon$ in $(0,\hat \varepsilon]$ and every $x\in \mathbb{B}(0,1)$, we have
\begin{equation}\label{equa:ref}
    1\ge\pi_i(x)\ge F_i^{\varepsilon}(x)\ge F_i^{\varepsilon_n}(x)\ge -1,
\end{equation}
therefore given $z^n$ such that 
\begin{equation}\label{equa:zeta_eps2}
F^{\varepsilon_n}(x)= \pi_i(z_n)+\frac{1}{\varepsilon_n}\rho(z_n,x),
\end{equation}
up to a subsequence, we have that $z_n\to \hat z$ and, by equation \eqref{equa:ref}, $\hat z$ realizes the minimum for $F^\varepsilon_i(x)$, thus we have $F_i^{\varepsilon_n}(x)\to F_i^{\varepsilon}(x).$
Therefore, for a generic sequence $\varepsilon_n\to\varepsilon$ in $(0,\hat \varepsilon]$ and every $x\in \mathbb{B}(0,1)$ fixed, up to subsequence we can assume $\varepsilon_n\searrow\varepsilon$ or $\varepsilon_n\nearrow\varepsilon$, hence we have that $F_i^{\varepsilon_n}(x)\to F_i^{\varepsilon}(x).$
In general, given $\varepsilon_n\to\varepsilon$ in $(0,\hat \varepsilon]$ and $x_n\to x$ in $\mathbb{B}(0,1)$, consider $z_n$ satisfying equation \eqref{equa:zeta_eps2} as before and observe that 
\begin{align}
    |F_i^{\varepsilon_n}(x_n)-F_i^{\varepsilon}(x)|&\le|F_i^{\varepsilon_n}(x_n)-F_i^{\varepsilon_n}(x)|+|F_i^{\varepsilon_n}(x)-F_i^{\varepsilon}(x)|\notag\\&\le\frac{\rho(x_n,x)}{\varepsilon_n}+o(1)=o(1)\qquad\text{for } n\to+\infty.
\end{align}
Finally, consider the last case when $\varepsilon_n\to 0$ and $x_n \to x$ in $\mathbb{B}(0,1)$, then 
\begin{equation}
    |F_i^{\varepsilon_n}(x_n)-\pi_i(x)|\le|F_i^{\varepsilon_n}(x_n)-\pi_i(x_n)|+|\pi_i(x_n)-\pi_i(x)|\le o(1),
\end{equation}
because of uniform convergence.\\
Now we consider the topological degree of the function $F^{\hat\varepsilon}$ with respect to the set $\mathbb{B}(0,1)$ and any point of $B(0,\frac{1}{2})$, the open ball of radius $\frac{1}{2}$. We recall that the map  $F^{\hat\varepsilon}$ is homotopy equivalent to the map id, and observe that equation \eqref{equa:bound} implies that for any $y\in B(0,\frac{1}{2})$, we have  $y\notin F^\varepsilon(\mathbb B(0,1)\setminus B(0,1))$. Therefore, we can apply homotopy invariance obtaining that
\begin{equation}
    1=\operatorname{deg}(\operatorname{id},\mathbb{B}(0,1), y)=\operatorname{deg}(F^{\hat \varepsilon},\mathbb{B}(0,1), y)
\end{equation} for every $y\in B(0,\frac{1}{2})$, hence, by Proposition \ref{prop degree},  it follows $B(0,\frac{1}{2})\subseteq F^{\hat \varepsilon}(\mathbb B(0,1))$ and this proves our claim.\\
Since $F^{\hat \varepsilon}(\mathbb B(0,1))$ contains a non-empty open set and $F^{\hat \varepsilon}$ is Lipschitz, we get
\begin{equation}
    \mathcal{H}^n_{\textrm{d}_\textrm{eucl}}\bigl(F^{\hat \varepsilon}(\mathbb{B}(0,1))\bigr)\leq C_\varepsilon^n\mathcal{H}^n_{\rho}(\mathbb{B}(0,1))=0,
\end{equation}
which is a contradiction since the n-dimensional Hausdorff measure on $\mathbb R^n$ with the Euclidean distance gives positive measure to not empty open sets.
\end{proof}
\begin{remark}
The same proof of Theorem \ref{main th} can be adapted to prove that any nonempty open set $A$ is such that $\mathcal{H}^n_\rho(A)>0$.
\end{remark}
\begin{remark}\label{pathologies}
Removing the assumption that $\rho$ induces the Euclidean topology, counterexamples show that $\mathcal{H}^n_{\rho}\left(\mathbb{R}^n\right)$ might vanish. Consider, for instance, the metric space $(\mathcal{C},\textrm{d})$, where  $\mathcal{C}\subset\mathbb{R}$ is the Cantor set and $\textrm{d}$ denotes the usual one-dimensional Euclidean distance. Having $\mathcal{C}$ the cardinality of the continuum, there exist bijections $g_n:\mathcal{C}\to\mathbb{R}^n$. Then, define on $\mathbb{R}^n$ the metric $\rho(x,y)=\textrm{d}(g_n^{-1}(x),g_n^{-1}(y))$.\\
Given any collection $(A_i)_{i\in\mathbb{N}}$ that covers $\mathcal{C}$, follows that $(g_n(A_i))_{i\in\mathbb{N}}$ covers $\mathbb{R}^n$ and $\operatorname{diam}(A_i)=\operatorname{diam}(g_n(A_i))\, \forall i \in \mathbb{N}$. Clearly, also the opposite direction applies. Therefore, we have
\begin{equation}
   \mathcal{H}^n_{\rho}\left(\mathbb{R}^n\right)= \mathcal{H}^n_{d}\left(\mathcal{C}\right)=0
\end{equation}
that shows a counterexample.
\end{remark}
\begin{remark}
Note that, under previous assumptions on $\rho$, it is not true in general that $\operatorname{dim}_H^{\rho}(\mathbb{R}^n)=n$. In fact, choosing $\rho(x,y)=\textrm{d}_{\textrm{eucl}}(x,y)^{1/2}$, the distance $\rho$ induces the Euclidean topology, but in this case
\begin{align*}
\mathcal{H}^s_{\textrm{d}_\textrm{eucl}}\left(A\right)=  \mathcal{H}^{2s}_{\rho}\left(A\right)
\end{align*}
for all $A\subseteq \mathbb{R}^n,\ s\geq 0$, see for example \cite{federer2014geometric}. For this reason we get that $\operatorname{dim}_H^{\rho}(\mathbb{R}^n)=2n$.
\end{remark}
\section*{Acknowledgements}
We would like to thank professor David H. Fremlin for proposing this problem and professor Stefano Bianchini for the inspiring and useful discussions we had about it.
\section*{Conflict of interest}
Authors state no conflict of interest.
  \bibliography{demo}

\begin{thebibliography}{1}

\bibitem{dinca2021brouwer}
G.~Dinca and J.~Mawhin.
\newblock {\em Brouwer Degree}.
\newblock Springer, 2021.

\bibitem{federer2014geometric}
H.~Federer.
\newblock {\em Geometric measure theory}.
\newblock Springer, 2014.

\bibitem{fremlin2000measure}
D.~H. Fremlin.
\newblock {\em Measure theory}, volume~4.
\newblock Torres Fremlin, 2000.

\bibitem{jost2005postmodern}
J.~Jost.
\newblock {\em Postmodern analysis}.
\newblock Springer Science \& Business Media, 2005.

\end{thebibliography}
  \bibliographystyle{abbrv}
\end{document}